\documentclass[10pt]{article}
 \usepackage{latexsym}
 \usepackage{amsbsy}
 \usepackage{amssymb}
 \usepackage{amsmath}
 \usepackage[latin1]{inputenc}
 \usepackage{graphicx}
 \def\dj{d\kern-0.4em\char"16\kern-0.1em}
 \def\Dj{\mbox{raise0.3ex\hbox{-}\kern-0.4em D}}

 \newcommand{\cc}{{\cal C}_{c}}
 \newcommand{\cjedan}{{\cal C}^{1}}

 \newcommand{\cbesk}{{\cal C}^{\infty}}
 
 \newcommand{\rp}{{\mathbb R}^p}
 \newcommand{\rqq}{{\mathbb R}^q}
 \newcommand{\rn}{{\mathbb R}^n}

 \newcommand{\test}{{\cal D}}
 \newcommand{\dprim}{{\cal D}'}

 \newcommand{\pd}{\partial}

 \newcommand{\vp}{{\bf v}}
 \newcommand{\vpw}{{\bf w}}

 \newcommand{\vpxu}{\vp=\sum_{i=1}^p\xi^i(x,u)\frac{\pd}{\pd x_i}+
 \sum_{\alpha=1}^q\phi_{\alpha}(x,u)\frac{\pd}{\pd u^{\alpha}}}

 \newcommand{\Un}{U^{(n)}}
 \newcommand{\un}{u^{(n)}}
 \newcommand{\prn}{\mbox{\rm pr}^{(n)}}
 
 \newcommand{\prjedan}{\mbox{\rm pr}^{(1)}}

 \newcommand{\eps}{\varepsilon}

 \newcommand{\R}{\mathbb R}
 \newcommand{\N}{\mathbb N}

 \newcommand{\cR}{{\cal R}}

 \newcommand{\cG}{{\cal G}}

 \newcommand{\be}{\begin{equation}}
 \newcommand{\ee}{\end{equation}}
 \newcommand{\bea}{\begin{eqnarray}}
 \newcommand{\eea}{\end{eqnarray}}
 \newcommand{\beast}{\begin{eqnarray*}}
 \newcommand{\eeast}{\end{eqnarray*}}
 \newcommand{\bary}{\begin{array}}
 \newcommand{\eary}{\end{array}}

 \newcommand{\pr}{{\bf Proof. }}
 \newcommand{\ep}{\hspace*{\fill}$\Box$}

 \newtheorem{thm}{\hspace*{-1.1mm}}[section]
 \newcommand{\bt}{\begin{thm} {\bf Theorem. }}
 \newcommand{\et}{\end{thm}}
 \newcommand{\bp}{\begin{thm} {\bf Proposition. }}
 \newcommand{\bc}{\begin{thm} {\bf Corollary. }}
 \newcommand{\blem}{\begin{thm} {\bf Lemma. }}
 \newcommand{\bex}{\begin{thm} {\bf Example. }\rm}
 \newcommand{\bexs}{\begin{thm} {\bf Examples. }\rm}
 \newcommand{\bd}{\begin{thm} {\bf Definition. }}
 \newcommand{\br}{\begin{thm} {\bf Remark. }\rm}
 \newcommand{\brs}{\begin{thm} {\bf Remarks. }\rm}

 \newcommand{\sms}{\smallskip\\}

\begin{document}
 \input{amssym.def}
 \input{amssym}

 \title{Symmetries of Conservation Laws}
 \author{Sanja Konjik \footnote{Electronic mail: kinjoki@neobee.net}\\
        {\small Departmant of Agriculture Engineering, Faculty of Agriculture}\\
        {\small University of Novi Sad}\\
        {\small Trg Dositeja Obradovica 8, 21000 Novi Sad}\\
        {\small Serbia and Montenegro}\\
        }

 \date{\mbox{}}
 \maketitle

 \begin{abstract}

 The aim of this paper is to apply techniques of symmetry group
 analysis in solving two systems of conservation laws: a model of
 two strictly hyperbolic conservation laws and a zero pressure gas
 dynamics model, which both have no global solution, but whose solution
 consists of singular shock waves. We show that these shock waves
 are solutions in the sense of $1$-strong association. Also, we compute
 all projectable symmetry groups and show that they are $1$-strongly
 associated, hence transform existing solutions in the sense of $1$-strong
 association into other solutions.

 \end{abstract}

 The concept of classical symmetry groups offers a large number of
 possibilities in studying differential equations, in particular
 in constructing explicit solutions to linear and nonlinear
 differential equations or determining and classifying invariance
 properties \cite{Olver1,Olver2}. In various problems of
 mathematical physics the classical theory turns out to be
 insufficient, due to singular objects (like distributions or
 discontinuous nonlinearities) which can occur in the equation or
 equations with solutions in a weak sense, i.e. weak solutions
 (distributional, generalized or in the sense of association).
 Therefore, the methods of classical symmetry group analysis of
 differential equations have been extended to linear equations in
 the class of distributions \cite{Ber1,Ber2}, as well as to equations
 involving generalized functions \cite{Dj-P-K,GKOS,Kunz-Ober1,Kunz-Ober2,KunzPhD}.

 The aim of this paper is to apply techniques of symmetry group
 analysis in solving two systems of differential equations given
 in the form of conservation laws. The paper is divided into two parts.
 Section $1$ provides a brief overview of the basic definitions and theorems
 which are going to be used for studying conservation laws. We
 start by recalling some facts on symmetry group analysis, which are in
 detail carried out in \cite{Olver1} (see also \cite{Olver2}). Then we turn to
 symmetries in the generalized setting, precisely to associated ones. As
 we will see later, the reason for this lies in the fact that the solutions of
 the conservation laws we consider are, under certain assumptions, shock
 waves. Lack of space prevents us from also giving a short
 introduction to generalized functions. Therefore, for the notations and
 properties of the Colombeau algebra of generalized functions we recommend
 \cite{GKOS,KunzPhD} or \cite{Ober-Kunz}; in particular, definitions of generalized
 step-functions, splitted delta functions and $m-$ and $m'-$singular delta
 functions are provided in \cite{Ned}. We close the introductory part by
 a short overview of conservation laws. Based on
 \cite{Bres1,Bres2,Hoe,LeFloch,Smol} we fix
 notations and present the general solution of the Riemann problem.
 Motivated by \cite{Keyfitz} and \cite{Ned} we proceed in
 section $2$ by investigating two systems of conservation laws: a
 model of two strictly hyperbolic conservation laws which is
 genuinely nonlinear but for which the Riemann problem has no
 global solution and a zero pressure gas dynamics model which is
 linearly degenerative but for which the Riemann problem also does
 not have global solutions. In both cases singular solutions appear,
 called singular shock waves. We prove that these solutions are solutions
 in the sense of $1$-strong association. After computing all projectable symmetry
 groups of these systems we show that they are $1$-strongly associated, hence
 transform existing solutions (given in \cite{Keyfitz} and \cite{Ned}) into
 other solutions.

\section{Introduction}

\subsection{Symmetry Groups of Differential Equations}

 Let $S$ be a system of differential equations:
 $$
 \Delta_{\nu} (x, u)=0, \quad\quad 1\leq \nu \leq l.
 $$
 Denote by
 $X=\rp$ and $U=\rqq$ the spaces of independent and dependent variables
 with coordinates $x=(x_1, x_2, \dots, x_p)$ and $u=(u^1, u^2, \dots, u^q)$
 respectively. Also, denote by $M$ an open subset of $X \times U$.
 Identify a function $u=f(x)$ with its graph
 $$
 \Gamma_f=\{(x, f(x)): x\in \Omega\} \subset X\times U,
 $$
 where $\Omega\subset X$ is the domain of $f$. Let $G$ be a local group
 of transformations acting on $M$. The transform of
 $\Gamma_f$ by $g\in G$ is defined by
 $$
 g\cdot\Gamma_f = \{(\tilde{x}, \tilde{u})=g\cdot(x,u): (x,u)\in
 \Gamma_f\}.
 $$
 In local coordinates this action is given by
 $$
 g\cdot f = \big(\Phi_g \circ (id_X \times f) \big)\circ\big(
 \Xi_g \circ (id_X \times f) \big)^{-1},
 $$
 where $\Xi_g$ and $\Phi_g$ are smooth function on $M$, and $id_X$ is the
 identity mapping on $X$. Supposing that $\Xi_g$
 does not depend on the dependent variables we get projectable action of
 $g$ on $f$, i.e.
 \be \label{ftilda projektabilna}
 g\cdot f = \big(\Phi_g \circ (id_X \times f)
 \big)\circ \Xi_g^{-1}.
 \ee

 \bd
 The symmetry group of the system $S$ is a local
 transformation group $G$ acting on the space of independent and
 dependent variables with the property that whenever $u=f(x)$
 is a solution of the system and $g\cdot f$ is defined, $g\in G$,
 then $u=g\cdot f$ is also a solution of the system.
 \et

 The $n$-th prolonged or $n$-jet space $X\times \Un$ is a space
 which represents all independent
 variables, dependent variables and all different partial derivatives
 of dependent variables up to the order $n$. For the construction
 of the $n$-th prolonged space we refer to \cite{Olver1}. Write
 $M^{(n)}$ for a subset of $n$-jet space $X\times \Un$.
 An arbitrary point in $\Un$ will be denoted by $u^{(n)}$ and its
 components by $u_J^{\alpha}$, where
 $1\leq \alpha \leq q$, while $J$ runs over the set of all unordered
 multi-indices $J=(j_1, \dots, j_k)$, $1\leq j_k \leq p$,
 $0\leq k\leq n$.

 The $n$-th prolongation of a function $f: X\to U$,
 denoted by $pr^{(n)} f$, is a function from $X$ to
 $\Un$, which maps $x$ into $ (\pd_J f^{\alpha}(x))_{\alpha, J}$,
 $1 \leq \alpha \leq q$, $0\leq |J| \leq n$.

 The $n$-th prolongation of a group $G$ which acts on
 $M\subset X\times U$, $\prn G$, is again a local group
 of transformations which acts on $M^{(n)}$ such that
 it transforms the derivatives of a smooth function $u=f(x)$ into
 the corresponding derivatives of the transformed function
 $\widetilde{u}=\widetilde{f}(\widetilde{x})$. For the precise
 definition see \cite{Olver1}.

 The $n$-th prolongation of a vector field $\vp$ on $M$, $\prn \vp$,
 is a vector field on the $n$-jet space $M^{(n)}$ with the following
 property:
 $$
 \prn \vp|_{(x, \un)} =
 \frac{d}{d\eta}{\Big\vert}_{_{\eta=0}} \prn (\exp(\eta\vp))(x, \un),
 $$
 where $\exp(\eta \vp)$ is the corresponding local one-parameter
 group generated by $\vp$.
 If
 $$
 \vpxu
 $$
 then we calculate the $n$-th prolongation of $\vp$ using the formula:
 \be \label{decembar2003}
 \prn \vp = \vp + \sum_{\alpha =1}^q \sum_{J}\phi_{\alpha}^J(x, \un)
 \frac{\pd}{\pd u^{\alpha}_J},
 \ee
 where the coefficients $\phi_{\alpha}^J(x, \un)$ are given by
 \be \label{prol formula za phi}
 \phi_{\alpha}^J(x, \un)=D_J \Bigg(
 \phi_{\alpha} - \sum_{i=1}^p \xi^i u_i^{\alpha}  \Bigg) +
 \sum_{i=1}^p \xi^i u_{J,i}^{\alpha},
 \ee
 $\displaystyle u_i^{\alpha} = \pd u^{\alpha}/{\pd x^i}$,
 $\displaystyle u_{J,i}^{\alpha} = \pd u_J^{\alpha}/{\pd x^i}$
 and $D_J$ denotes a total differential.

 Then the infinitesimal criterion for a system of differential
 equations reads:

 \bt \label{thm infinitesimalni kriterijum}
 Let
 \be \label{sistem pkzkp}
 \Delta_\nu(x, \un)=0, \quad\quad \nu=1,\dots,l
 \ee
 be a system of differential equations of a maximal rank
 (meaning that the corresponding
 Jacobian matrix
 $ \mbox{J}_{\Delta}(x, \un)=(\pd \Delta_{\nu}/\pd x^i,
 \pd \Delta_{\nu}/\pd u^{\alpha}_J)$ is of rank $l$ on the set
 of all solutions of $S$, $S_{\Delta}$).
 If $G$ is a local transformation group acting on
 $M\subset X\times U$ and
 \be \label{infinitesimalni kriterijum}
 \prn \vp (\Delta_{\nu}(x, \un))=0, \quad \nu=1, \dots,l,
 \quad \mbox{ whenever } \quad \Delta(x, \un)=0,
 \ee
 for every infinitesimal generator $\vp$ of $G$
 then $G$ is a symmetry group of (\ref{sistem pkzkp}).
 \et

 The condition (\ref{infinitesimalni kriterijum}) from this theorem
 will also be necessary if we additionaly suppose that the 
 system (\ref{sistem pkzkp}) is locally solvable, i.e. at each
 point $(x_0, \un_0)\in S_{\Delta}$ there exists a smooth
 solution $u=f(x)$ of the system, defined in a neighborhood of
 $x_0$, which has the prescribed ``initial conditions'' $\un_0 
 = \prn f(x_0)$. (We say that a system of differential equations
 is nondegenerate if at every point of the solution set it is both
 locally solvable and of maximal rank.) 

 For later use, we mention here a result which is a
 consequence of the maximal rank condition (imposed on the system
 (\ref{sistem pkzkp}) in the above theorem). Namely, under the
 conditions of Theorem \ref{thm infinitesimalni kriterijum}, the
 infinitesimal criterion (\ref{infinitesimalni kriterijum}) can be
 replaced by the equivalent condition
 \be \label{ekv infkr}
 \prn \vp (\Delta_{\nu}(x, \un))=\sum_{\mu=1}^{l} Q_{\nu\mu}
 \Delta_{\mu}(x, \un), \quad \nu=1, \dots,l,
 \ee
 for functions $Q_{\nu\mu}$, $\mu,\nu=1,\dots,l$ to be determined. 

 We finish this short introduction into symmetry groups of
 differential equations by a description of a procedure for calculating
 symmetry groups of a given system $S$. The procedure consists of
 the following steps:

 \begin{itemize}
 \item[(1)] Write the vector field (i.e. infinitesimal generator)
 in the most general form:
 $$
 \vp (x,u) = \sum\limits^p_{i=1}\xi_i(x,u)\frac{\pd}{\pd x_i} +
 \sum\limits^q_{\alpha =1} \phi_\alpha (x,u) \frac{\pd}{\pd
 u^\alpha},
 $$
 where $\xi_i$ and $\phi_\alpha$ are functions which should be calculated.
 \item[(2)] According to (\ref{decembar2003}) and (\ref{prol formula za phi})
 calculate the corresponding prolongation of $\vp$.
 \item[(3)] Then apply the infinitesimal criterion
 (\ref{infinitesimalni kriterijum}) and equate
 $\prn \vp (\Delta_\nu(x, \un))$ with zero. Since those equations
 must hold on $S_\Delta$, eliminate the dependence of derivatives of $u$
 by the equations from the system. After that we have the equations which
 have to be satisfied with respect to $x$, $u$ and the remaining partial
 derivatives of $u$.
 \item[(4)] After solving these equations we obtain a certain number
 of partial differential equations for $\xi_i$ and $\phi_\alpha$.
 \item[(5)] Compute the $\xi_i$ and $\phi_\alpha$ from them,
 thereby computing vector fields $\vp$ which generate
 a Lie algebra of infinitesimal symmetries.
 \item[(6)] At the end find the corresponding one-parameter symmetry
 groups as the flows of the infinitesimal generators calculated in the
 previous step.
 \end{itemize}

\subsection{Symmetry Groups of Weak Solutions}

 Next, we look for the symmetries which transform weak solutions
 of the system of PDEs
 \be \label{ssiisstteemm}
 \Delta_{\nu} (x, \un)=0, \quad\quad 1\leq \nu \leq l,
 \ee
 into other weak solutions, mainly associated solutions to
 (\ref{ssiisstteemm}) into other associated solutions to
 (\ref{ssiisstteemm}) (hence the system (\ref{ssiisstteemm})
 should be replaced by $\Delta_{\nu} (x, \un)\approx 0$).
 Such a symmetry group is called symmetry group in the sense of
 association or associated symmetry group for short.
 The symmetry groups we are interested in are projectable in order
 to avoid the problem of inverting Colombeau functions.
 Thus unless explicitly stated
 otherwise, all symmetry groups are assumed to be projectable.
 Beside this, we need some more assumptions on $G$ (cf.\ \cite{Dj-P-K}).
 We suppose that a local transformation group $G$ is
 slowly increasing, uniformly for $x$ in compact sets, and 
 analogously for the mapping $ \un \mapsto \Delta(x, \un)$.

 Look at the system
 \begin{equation}\label{asssystem}
 \Delta_\nu(x,\un) \approx 0,  \qquad \nu = 1,\dots,l.
 \end{equation}
 We recall the following definitions from \cite{Dj-P-K}:
 \bd \label{asocirano resenje}
 $u=(u^1,...,u^{q}) \in {\cal{G}}(\Omega )^{q}$ is a solution of
 (\ref{asssystem}) and also associated solution to
 (\ref{ssiisstteemm}) if $u=(u^1,...,u^{q})$ has a representative
 $[(u^1_{\eps},...,u^q_{\eps})_\eps] \in {\cal{E}}_{M}(\Omega )^q$
 such that for each test function $\varphi \in \test(\Omega)$
 \begin{equation} \label{assdef}
 \int \Delta _{\nu} (x,\prn u_\eps(x))\varphi (x)\, dx \rightarrow 0
 \quad \mbox{ as }\;\; \eps \rightarrow 0, \;\; 1\le \nu \le l.
 \end{equation}
 \et

 The set of all associated solutions to (\ref{ssiisstteemm})
 is denoted by ${\cal{A}} _{\Delta }$ and moreover, the set of all
 $u \in (\cG_{\infty})^q$ which satisfies (\ref{assdef}) with
 ${\cal{A}}{\cal{B}}_{\Delta }$. The symmetry group $G$ of
 (\ref{ssiisstteemm}) is called ${\cal{A}}$-symmetry group
 if for every $u\in {\cal A}$ and every $g_\eta \in G$ it
 follows that $g_\eta u \in {\cal A}$, whenever $g_\eta
 u$ is defined. Beside solution in the sense of association
 we can also define a solution in the sense of strong association.

 \bd \label{def strass}
 Let $k\in \N_0$.
 $u = (u^{1},...,u^{q}) \in {\cal{A}} _{\Delta }(\Omega)$
 (resp. $u \in {\cal{A}}{\cal{B}}_{\Delta }(\Omega))$ is called
 $k$-strongly associated or $\stackrel{k}{\approx}$-associated
 solution to the system (\ref{ssiisstteemm})
 if there exists a representative
 $[(u^1_{\epsilon },...,u^q_{\epsilon })_\eps] \in
 {\cal{E}}_{M}(\Omega )^{q}$ such that for each
 $B \subseteq C^{\infty }_c (\Omega )$ which is bounded in
 ${\cal{C}}_c^k(\Omega)$ we have
 $$
 \lim _{\epsilon \rightarrow 0} \sup _{\varphi \in B}|\int \Delta_{\nu}
 (x,\prn u_\eps(x))\varphi (x)\, dx| = 0, \qquad 1\le \nu\le l.
 $$
 \et

 The space of all $k$-strongly associated solutions to
 (\ref{ssiisstteemm}) is denoted by ${\cal AS}^k_\Delta$.
 Also,
 ${\cal ABS}_\Delta^k:={\cal AS}^k_\Delta\cap \cG_\infty$.
 The main role in the calculation of associated, resp.
 $k$-strongly associated symmetry groups is played by the theorem
 which is based on the factorization property of the system,
 derived in \cite{GKOS}. In matrix form this property
 is given by

 \begin{equation} \label{faktorizacija matricna}
 \Delta(\Xi_\eta(x,u(x)), \prn(g_\eta u)(\Xi_\eta(x,u(x))))
 \!=\!Q(\eta,x,\prn u(x))\Delta(x,\prn u(x)),
 \end{equation}

 where $Q:{\cal W} \to \R^{l^2}$ and $\cal W$ is an open
 subset of $(-\eta_0, \eta_0)\times M^{(n)}$ with
 $\{0\}\times M^{(n)} \subseteq \cal W$.

 \bt \label{thm Q}
 Let $G$ be a slowly increasing symmetry group of the system
 (\ref{ssiisstteemm}) which admit a global factorization of the form
 (\ref{faktorizacija matricna}). Then
 \begin{itemize}
 \item[(i)] if $Q$ depends only on  $\eta$, $x$ and $u$ then
 $G$ is also an ${\cal ABS}_\Delta$-symmetry group of
 (\ref{ssiisstteemm});
 \item[(ii)] if $Q$ depends only on $\eta$ and $x$ then
 $G$ is also an ${\cal AS}_\Delta^k$-symmetry group of
 (\ref{ssiisstteemm}), for each $k>0$. Moreover,
 $G$ is in this case an associated symmetry group of
 (\ref{ssiisstteemm}) as well.
 \end{itemize}
 \et

 \subsection{Systems of Conservation Laws}

 We look at a system of conservation laws in one space dimension:
 \bea \label{systemcl}
 u^1_t + (f_1(u^1, \dots, u^n))_x & = & 0 \nonumber \\
 \vdots                          &   & \nonumber \\
 u^n_t + (f_n(u^1, \dots, u^n))_x & = & 0, \nonumber
 \eea
 or written in a shorter (matrix) form:
 \be \label{cl}
 u_t + (f(u))_x=0,
 \ee
 where $t>0$, $x\in \R$, $u=(u^1, \dots , u^n)$
 is the conserved density and
 $f(u)=(f_1, \dots , f_n)$ is the flux. Differentiating (\ref{cl})
 we obtain a quasilinear system
 \be \label{smoothsystemcl}
 u_t + A(u)u_x =0,
 \ee
 where $A(u) = Df(u)$ is the Jacobian matrix of $f$.
 The systems (\ref{cl}) and (\ref{smoothsystemcl})
 are equivalent for all smooth solutions $u$. Otherwise, if
 $u$ has a jump, the left hand side of
 (\ref{smoothsystemcl}) contains a product of a discontinuous
 function with a distributional derivative, while (\ref{cl})
 is still well defined in the distributional sense.

 The eigenvalues of the matrix $A(u)$ determine the system of
 conservation laws in the following way:

 \bd
 The system of conservation laws is hyperbolic, resp. strictly
 hyperbolic, if all eigenvalues of the matrix $A(u)$ are
 real, resp. real and different.
 \et

 Suppose that the system (\ref{cl}) is strictly hyperbolic
 and denote by $\lambda_1(u),
 \dots,$ $\lambda_n(u)$ the eigenvalues of $A(u)$ with
 $\lambda_1(u) < \dots < \lambda_n(u)$. Next, denote by
 $l_1, \dots, l_n$ and $r_1, \dots, r_n$ the corresponding
 left and right eigenvectors. The eigenvalue $\lambda_i$ of $A$
 is also called the $i$-th characteristic speed and the pair
 $(\lambda_i, r_i)$ the $i$-th characteristic field of (\ref{cl}).

 \bd  \label{gen nonlin and lin deg}
 The $i$-th characteristic field of the system (\ref{cl})
 is called genuinely nonlinear if
 $$
 D\lambda _i(u) \cdot r_i(u)\not= 0, \qquad \forall u.
 $$
 The $i$-th characteristic field is called
 linearly degenerate if
 $$
 D\lambda _i(u) \cdot r_i(u)= 0, \qquad \forall u.
 $$
 \et

 If the solution of (\ref{cl}) is a piecewise smooth function
 $u=u(t,x)$ having a discontinuity across a line $x=\gamma (t)$ with
 $u^{\pm}=\lim_{x\to \gamma(t)^{\pm}} u(t,x)$, then it satisfies
 (\ref{smoothsystemcl}) outside the $\gamma$, while along the
 line of discontinuities the Rankine-Hugoniot conditions holds:
 $$
 (u^+ - u^-)\dot{\gamma}=f(u^+) - f(u^-).
 $$

 In order to have a unique solution we must require some
 additional conditions, known as entropy conditions.
 One of the most useful is the Lax condition, which says that a
 shock connecting the states $u^-$ and $u^+$, travelling with speed
 $\dot{\gamma}=\lambda_i(u^-, u^+)$ ($\lambda_i(u^-, u^+)$ is an
 eigenvalue of the averaged matrix $A(u^-, u^+)$, cf.\ \cite{Bres1})
 is admissible if
 \be \label{Laxov uslov}
 \lambda_i(u^-) \geq \lambda_i(u^-, u^+) \geq
 \lambda_i(u^+).
 \ee

 Now we define two types of curves: fix a state $u_0\in \rn$.
 Let $r_i(u)$ be the $i$-th eigenvector of $A(u)$. The
 $i$-th rarefaction curve through $u_0$ is the integral
 curve of the vector field $r_i$ and is denoted by
 $\sigma \mapsto R_i(\sigma)(u_0)$. The
 $i$-th shock curve through $u_0$ is the curve of states $u$
 which can be connected to the right of $u_0$ by an $i$-shock,
 satisfying the Rankine-Hugoniot conditions. It is denoted by
 $ \sigma \mapsto S_i(\sigma)(u_0)$.
 The $i$-th rarefaction and shock curve are tangent to the
 $r_i(u)$ at $u_0$.

 Next we study the Riemann problem

 \begin{equation} \begin{split} \label{rimanov problem}
 u_t + f(u)_x &= 0 \\
 u(0,x)       &= \left\{\begin{array}{ll}
                          u^-, & \; x<0,\\
                          u^+, & \; x>0.
                          \end{array}
                          \right.
 \end{split} \end{equation}

 Under the assumption that the system is strictly hyperbolic
 with smooth coefficients and each $i$-th characteristic field
 is either genuinely nonlinear or linearly degenerate, there exist
 three special cases:

 \begin{itemize}
 \item[(1)] Centered rarefaction waves: the $i$-th characteristic
 field is genuinely nonlinear and $u^+$ lies on the positive
 $i$-rarefaction curve through $u^-$, i.e. $u^+=R_i(\sigma)(u^-)$ for
 some $\sigma>0$. Then the solution of (\ref{rimanov problem})
 is the centered rarefaction wave:
 $$
 u(t,x) = \left\{\begin{array}{ll}
                   u^-, & \quad\quad x<t\lambda_i(u^-),\\
                   R_i(s)(u^-), & \quad\quad
                   x=t\lambda_i(s), \;\;\;\;s\in[0, \sigma]\\
                   u^+, & \quad\quad x>t\lambda_i(u^+).
                          \end{array}
                          \right.
 $$

 \item[(2)] Shocks: again the $i$-th characteristic
 field is genuinely nonlinear, but this time the state
 $u^+$ lies on the $i$-th shock curve through $u^-$,
 i.e. $u^+=S_i(\sigma)(u^-)$. Denote the Rankine-Hugoniot
 speed of the shock $\lambda_i(u^-, u^+)$ by $\lambda$. Then
 the solution of (\ref{rimanov problem}) is the shock
 \be \label{res s}
 u(t,x) = \left\{\begin{array}{ll}
                   u^-, & \quad x<\lambda t,\\
                   u^+, & \quad x>\lambda t.
                          \end{array}
                          \right.
 \ee

 \item[(3)] Contact discontinuities: the $i$-th characteristic
 field is linearly degenerate and $u^+$ lies on the $i$-th
 rarefaction curve through $u^-$, i.e. $u^+=R_i(\sigma)(u^-)$
 for some $\sigma$. Then the function (\ref{res s}) is again a solution,
 but this time called the contact discontinuity.
 \end{itemize}

 The parameter $\sigma$, for which $u^+=R_i(\sigma)(u^-)$ or
 $u^+=S_i(\sigma)(u^-)$, is called the wave strength.

 Therefore, if $u^+$ lies on the rarefaction or shock curve the
 solution of the Riemann problem
 (\ref{rimanov problem}) is one of the elementary waves -
 a centered rarefaction, a shock or a contact discontinuity.
 Otherwise, for $u^+$ sufficiently close to $u^-$,
 the Riemann problem (\ref{rimanov
 problem}) can be decomposed in $n$ auxiliary Riemann
 problems, which can be solved by an elementary wave. Piecing
 together those solutions we obtain a solution of the initial
 Riemann problem (\ref{rimanov problem}).

\section{Symmetry Groups of the Systems of Conservation Laws}

 After we gave the brief overview of the notation and results
 from symmetry group analysis (classical and in the generalized
 setting) and conservation laws, we turn our attention to concrete
 systems of conservation laws. As we have just seen, a system of
 conservation laws is a system of first order partial differential
 equations. We introduced symmetry groups of system
 of differential equations as local transformation groups which
 act on the space of independent and dependent variables,
 transforming the solution of the system to other solutions.
 Also, we defined associated and $k$-strongly associated
 symmetry groups. The aim of this section is to verify the results
 given in the introduction in two examples of conservation laws.

\subsection{A Model System of Two Strictly Hyperbolic Laws}
\label{section prvi sistem}

 The first system we consider is

 \begin{equation} \begin{split}\displaystyle \label{sistem1}
 u_t + (u^2 -v)_x & = 0 \\
 v_t + \Big(\displaystyle \frac{1}{3}u^3-u\Big)_x & = 0
 \end{split} \end{equation}

 with initial conditions:
 \be \label{pocetni uslovi}
 u(x,0)=\left\{\begin{array}{rl}
      u_{0}, & \; x<0 \\
      u_{1}, &\; x>0
      \end{array}\right. \qquad
 v(x,0)=\left\{\begin{array}{rl}
      v_{0}, &\; x<0 \\
      v_{1}, &\; x>0.
      \end{array}\right.
 \ee

 A motivation for studying this system arises from some physical
 models like a model for a nonlinear elastic system or a model
 for the evolution of ion-acoustic waves.

 We start by calculating symmetry groups of the system
 (\ref{sistem1}), using the procedure described in the
 introduction.

 (1) Since $v$ denotes one of the dependent variables,
 denote the infinitesimal generator by ${\bf w}$:
 $$
 \vpw = \xi(x,t,u,v)\pd_x + \tau(x,t,u,v)\pd_t +
 \phi(x,t,u,v)\pd_u + \psi(x,t,u,v)\pd_v.
 $$

 (2) From (\ref{decembar2003}) we get the first prolongation
 of this vector field
 $$
 \prjedan \vpw = \vpw + \phi^x\pd_{u_x} + \phi^t\pd_{u_t} +
 \psi^x\pd_{v_x} + \psi^t\pd_{v_t},
 $$
 and by (\ref{prol formula za phi}) we calculate
 \be \label{guzva}
 \begin{array}{c@{\,=\,}c@{\,+\,}c@{\,+\,}c@{\,-\,}c@{\,-\,}c@{\,-\,}
 c@{\,-\,}c@{\,-\,}c@{\,-\,}c} \!\!\!\!\!
 \phi^x & \phi_x & \phi_u
 u_x & \phi_v v_x & \xi_x u_x & \xi_u u_x^2
 & \xi_v u_x v_x & \tau_x u_t & \tau_u u_x u_t & \tau_v u_t v_x \\
 \!\!\!\!\! \phi^t & \phi_t & \phi_u u_t & \phi_v v_t & \xi_t u_x &
 \xi_u u_x u_t
 & \xi_v u_x v_t & \tau_t u_t & \tau_u u_t^2 & \tau_v u_t v_t \\
 \!\!\!\!\! \psi^x & \psi_x & \psi_u u_x & \psi_v v_x & \xi_x v_x &
 \xi_u u_x v_x
 & \xi_v v_x^2   & \tau_x v_t & \tau_u u_x v_t & \tau_v v_x v_t \\
 \!\!\!\!\! \psi^t & \psi_t & \psi_u u_t & \psi_v v_t & \xi_t v_x &
 \xi_u u_t v_x
 & \xi_v v_x v_t & \tau_t v_t & \tau_u u_t v_t & \tau_v v_t^2 \\
 \end{array}
 \ee

 (3) Now we have
 \beast
 & \Delta_1(x,t,u,v,u_x,v_x,u_t,v_t) = u_t + 2u u_x - v_x &\\
 & \Delta_2(x,t,u,v,u_x,v_x,u_t,v_t) = v_t + u^2u_x - u_x, &
 \eeast
 therefore we need to solve the system
 \beast
 & \mbox{\rm pr}^{(1)}{\bf w}(\Delta_1) = \phi^t + 2\phi u_x +
         2\phi^x u - \psi^x = 0 &\\
 & \mbox{\rm pr}^{(1)}{\bf w}(\Delta_2) = \psi^t + 2\phi u u_x +
         \phi^x u^2 - \phi^x = 0, &
 \eeast
 whenever $u_t=-2uu_x+v_x$ and $v_t=-u^2u_x + u_x$.
 Inserting (\ref{guzva}) in these equations and replacing
 $u_t$ by $-2uu_x+v_x$ and $v_t$ by $-u^2u_x + u_x$, whenever
 they appear, we arrive at
 \beast
 &&\phi_t + \phi_u (-2uu_x+v_x) + \phi_v (-u^2u_x + u_x) - \xi_t u_x
 - \xi_u u_x (-2uu_x+v_x) \\
 && - \xi_v u_x (-u^2u_x + u_x)
 -\tau_t (-2uu_x+v_x)
 - \tau_u (-2uu_x+v_x)^2  \\
 && - \tau_v (-2uu_x+v_x) (-u^2u_x + u_x)+
 2\phi u_x +2u [ \phi_x + \phi_u u_x +
 \phi_v v_x \\
 && - \xi_x u_x - \xi_u u_x^2
 - \xi_v u_x v_x -
 \tau_x (-2uu_x+v_x) - \tau_u u_x (-2uu_x+v_x)\\
 && - \tau_v (-2uu_x+v_x) v_x ] -
 [ \psi_x + \psi_u u_x + \psi_v v_x
 - \xi_x v_x - \xi_u u_x v_x  \\
 && - \xi_v v_x^2 - \tau_x (-u^2u_x +
 u_x)- \tau_u u_x (-u^2u_x + u_x) - \tau_v v_x (-u^2u_x + u_x)]\\
 &&= 0\\
 \eeast
 \beast
 &&\psi_t + \psi_u (-2uu_x+v_x) + \psi_v (-u^2u_x + u_x) - \xi_t v_x -
 \xi_u (-2uu_x+v_x) v_x \\
 && - \xi_v v_x (-u^2u_x + u_x) - \tau_t (-u^2u_x + u_x)
 - \tau_u (-2uu_x+v_x) (-u^2u_x + u_x) \\
 && - \tau_v (-u^2u_x + u_x)^2 + 2\phi u u_x + (u^2-1)
 [\phi_x + \phi_u u_x + \phi_v v_x -
 \xi_x u_x - \xi_u u_x^2 \\
 &&- \xi_v u_x v_x -
 \tau_x (-2uu_x+v_x) - \tau_u u_x (-2uu_x+v_x)
 - \tau_v (-2uu_x+v_x) v_x ]\\
 &&=0.
 \eeast

 (4) Apparently, solving this system is quite complicated.
 Hence, we are going to look only for projectable symmetry groups.
 So, assume that $\xi$ and $\tau$ only depend on $x$ and $t$.
 Then we have
 \beast
 && \phi_t - 2\phi_uuu_x + \phi_uv_x - \phi_vu^2u_x \phi_vu_x -
 \xi_tu_x + 2\tau_tuu_x - \tau_tv_x + 2\phi u_x  \\
 && +\, 2\phi_xu + 2\phi_uuu_x + 2\phi_vuv_x - 2\xi_xuu_x + 4\tau_xu^2u_x -
 2\tau_xuv_x - \psi_x - \psi_uu_x  \\
 && -\, \psi_vv_x + \xi_xv_x -\tau_xu^2u_x + \tau_xu_x = 0 \\
 && \psi_t - 2\psi_uuu_x + \psi_uv_x - \psi_vu^2u_x + \psi_vu_x -
 \xi_tv_x + \tau_tu^2u_x - \tau_tu_x + 2\phi uu_x  \\
 && +\, \phi_xu^2 + \phi_uu^2u_x + \phi_vu^2 v_x - \xi_xu^2u_x + 2\tau_xu^3u_x
 - \tau_xu^2v_x - \phi_x - \phi_uu_x  \\
 && -\, \phi_vv_x + \xi_xu_x - 2\tau_xuu_x + \tau_xv_x = 0.
 \eeast
 These equations are in fact polynomials of free variables
 $x, t, u, v, u_x$ and $v_x$. The solution will be found
 by looking at their coefficients on the left and right hand
 side of the equations. Since the functions $\phi$, $\psi$
 and their derivatives depend on $x, t, u$ and $v$ we equate
 the coefficients of $1$, $u_x$ and $v_x$ to $0$. Then
 we arrive at the following equivalent system:
 \bea
 && \phi_t - \psi_x + 2u \phi_x = 0 \nonumber\\
 && 2\phi + \phi_v - \psi_u - \xi_t + \tau_x +
    u(2\tau_t - 2\xi_x) + u^2(3\tau_x-\phi_v) = 0 \nonumber\\
 && \phi_u - \psi_v + \xi_x - \tau_t +u(2\phi_v-2\tau_x)=0 \nonumber\\
 && \psi_t - \phi_x + u^2\phi_x = 0 \nonumber\\
 && \psi_v - \phi_u - \tau_t + \xi_x + u(2\phi - 2\psi_u - 2\tau_x)
    + u^2(\phi_u - \psi_v - \xi_x + \tau_t) + 2u^3\tau_x =0 \nonumber\\
 && \psi_u - \phi_v - \xi_t + \tau_x + u^2(\phi_v-\tau_t)=0 \nonumber
 \eea

 (5) The general solution of this system is:
 \bea
 && \xi(x,t) = c_1 x + c_3 t + c_4 \nonumber \\
 && \tau(x,t) = c_1 t + c_2  \nonumber \\
 && \phi(x,t,u,v) = c_3  \nonumber \\
 && \psi(x,t,u,v) = c_3 u + c_5 \nonumber
 \eea
 $c_1 - c_5$ are arbitrary constants. The linearly independent
 infinitesimal generators of projectable symmetry groups are:
 \beast
 && {\bf w}_1 = x \partial_x + t \partial_t \\
 && {\bf w}_2 = \partial_t \\
 && {\bf w}_3 = t \partial_x + \partial_u + u\partial_v \\
 && {\bf w}_4 = \partial_x \\
 && {\bf w}_5 = \partial_v
 \eeast

 (6) It remains to compute the corresponding one-parameter
 symmetry groups. Hence, the one-parameter group $G_1$, generated
 by the vector field $\vpw_1$, is a solution to the system of ODEs
 \beast
 && \dot{x}(\eta) = x(\eta) \\
 && \dot{t}(\eta) = t(\eta),
 \eeast
 with initial data $x(0)=x$ and $t(0)=t$, i.e.
 \[
 g_\eta: (x,t,u,v) \to \left(e^{\eta}x,e^{\eta}t, u,v \right).
 \]
 Since $G_1$ is a symmetry group, (\ref{ftilda projektabilna})
 implies that if  $u$ and $v$ are solutions of (\ref{sistem1})
 so are the functions
 \beast
 && \widetilde{u}: (x,t) \to u\left(e^{-\eta}x,e^{-\eta}t\right) \\
 && \widetilde{v}: (x,t) \to v\left(e^{-\eta}x,e^{-\eta}t\right).
 \eeast
 We repeat the same procedure for the remaining symmetry groups
 and calculate that $G_2, G_4$ and $G_5$ are translations of $t$,
 $x$ and $v$, respectively. Finally, the action of $G_3$ is given by
 \[
 g_\eta: (x,t,u,v) \to \left(x+\eta t,t,
 u+\eta,v + \eta u +\frac{\eta^2}{2} \right),
 \]
 and the functions
 \beast
 && \widetilde{u}: (x,t) \to u(x-\eta t, t)+\eta\\
 && \widetilde{v}: (x,t) \to v(x-\eta t,t)+\eta u(x-\eta t,t) +
       \frac{3\eta^2}{2}
 \eeast
 are solutions of the system whenever $u$ and $v$ are.
 \sms

 Therefore, we calculated all projectable symmetry groups of the
 system (\ref{sistem1}) and all transformed solutions.
 As we saw, the calculation of non-projectable symmetry groups
 is rather complicated on the one hand, and on the other,
 as was mentioned in the introduction, it is enough to study
 projectable groups if the solution is in $\cG$, $\dprim$ or
 if it is a solution in the sense of association.
 We recall from \cite{Keyfitz} that the
 system (\ref{sistem1}) has a solution in the sense of association:
 the Jacobian matrix of (\ref{sistem1}) is
 $$
 A = Df = \left[\begin{array}{cc}
           2u & -1 \\
           u^2 -1 & 0
          \end{array} \right],
 $$
 the eigenvalues are $\lambda_1(u, v)= u-1$ and
 $\lambda_2(u, v)= u+1$, and the corresponding right
 eigenvectors are $r_1(u, v)=[1\quad u+1]^T$ and
 $r_2(u, v)=[1\quad u-1]^T$. Since
 $D\lambda_i(u, v)\cdot r_i(u, v) > 0$, $i=1,2$, it follows
 that both characteristic fields are genuinely nonlinear
 and hence the solution consists only of centered rarefaction
 waves and shocks.

 The rarefaction curves are calculated as the integral
 curves of the vector fields $r_1$ and $r_2$:
 \bea
 && R_1 = \{(u,v): v=\frac12 u^2 + u + c_1\} \nonumber \\
 && R_2 = \{(u,v): v=\frac12 u^2 - u + c_2\}, \nonumber
 \eea
 and the shocks are found from the Rankine-Hugoniot equations:
 \be \label{sistem1 sok krive}
 v-v_0 = (u-u_0)\Big(\frac{u+u_0}{2}\mp \sqrt{1-\frac{(u-u_0)^2}{12}}\Big),
 \quad \mbox{for}\;\; |u-u_0| \leq 12.
 \ee

 The corresponding shock speeds are
 $$
 \dot{\gamma}= u_0 + \frac{u-u_0}{2}\pm \sqrt{1-\frac{(u-u_0)^2}{12}}\ ,
 $$
 where the sign $-$ refers to the speed of $1$-shock,
 and $+$ to the speed of $2$-shock. The Riemann problem
 (\ref{sistem1})-(\ref{pocetni uslovi}) has a classical
 solution for each $(u, v)$ lying in the area bounded by
 \beast
  J(u,v)& =& \Big\{(u,v): v=\frac{1}{2} u^2 + u + \frac{9}{2} + v_0 -
 \frac{1}{2} u^2_0 - u_0 \land u\geq u_0-3 \Big\}\\
 J_1(u,v)& = & \Big\{(u,v): (u,v)\ \mbox{satisfies}\ (\ref{sistem1 sok krive})
 \wedge u \leq u_0 - 3\Big\} \\
 J_2(u,v)& = & \Big\{(u,v): v=\frac{1}{2} u^2 - u - \frac{9}{2} + v_0 -
 \frac{1}{2} u^2_0 + u_0 \wedge u\geq u_0-3\Big\}
 \eeast
 The remaining $(u,v)$ are in the exterior of this area,
 which we denote by $Q$ and which is divided by the curves
 \beast
 D(u,v)& =& \{(u,v): v=v_0+ u^2 +(1-u_0) u - u_0 \land u\leq u_0-3 \}\\
 E(u,v)& = &\{(u,v): v=v_0+(u-u_0)(u_0-1) \wedge u\leq u_0-3 \}
 \eeast
 into three open regions. In each of them the
 solution consists of a singular shock wave which is given by:

 \begin{equation} \begin{split}\label{prvi sistem kolombo resenje}
 U(x,t) & = G(x-ct)+s_{1}(t)(\alpha_{0}d^{-}(x-ct)+
 \alpha_{1}d^{+}(x-ct))\\
 V(x,t) & = H(x-ct)++s_{2}(t)(\beta_{0}D^{-}(x-ct)+
 \beta_{1}D^{+}(x-ct)),
 \end{split} \end{equation}

 where $G(x-ct)$ and $H(x-ct)$ are generalized step functions
 (cf.\ \cite{Ned}, Def.\ 1(a)),
 $D(x-ct)=\beta_0 D^-(x-ct) + \beta_1 D^+(x-ct)$ is an
 S$\delta$-function with value $(\beta_0, \beta_1)$, $\beta_0
 + \beta_1 = 1$ (cf.\ \cite{Ned}, Def.\ 1(b)),
 $d(x-ct)=\alpha_0 d^-(x-ct)+\alpha_1 d^+(x-ct)$ is an
 $3'$SD-function with value $(\alpha_0, \alpha_1)$
 (cf.\ \cite{Ned}, Def.\ 3 with Ex.\ (ii)), such that
 $D(x-ct)$ and $d(x-ct)$ are compatible and
 \bea
 &&-c[G]+[G^2]-[H]=0, \label{prvi sistem uslov1}\\
 &&s_{2}(t)=s_{1}^2(t)(\alpha_0^2+\alpha_1^2),\label{prvi sistem uslov2}\\
 &&s_{2}(t)=\sigma_{1}t,\; \sigma_{1}=c[H]-\frac13[G^3]+[G],\;\sigma_{1}>0,
     \label{prvi sistem uslov3}\\
 &&cs_2(t)=s_1^2(t)(\alpha_0u_0+\alpha_1u_1).\label{prvi sistem uslov4}
 \eea
 The function $s_2(t)$ is called the strength of the singular
 shock wave and is the most important part of the solution
 which has to be uniquely determined.
 $\alpha_0$ and $\alpha_1$ can be chosen such that
 $\alpha_0^2+\alpha_1^2=1$, hence the condition (\ref{prvi sistem
 uslov2}) becomes
 $$
 s_2(t)=s_1^2(t),
 $$
 and the condition (\ref{prvi sistem uslov4})
 $$
 \alpha_0^2u_0+\alpha_1^2u_1=c.
 $$

 We are going to show that this is a solution in the sense of
 $1$-strong association.

 \bt \label{thm sistem1 1strogo asoc}
 The solution (\ref{prvi sistem kolombo resenje})
 of the system (\ref{sistem1}) is a $1$-strongly associated solution
 to (\ref{sistem1}).
 \et

 \pr
 In order to show this we use Definition
 \ref{def strass}. Let $B$ be a bounded subset of $\cc^1(\R\times
 [0,\infty))$. This means that there exists $K\subset\subset
 \R\times[0,\infty)$, with the property $\mbox{supp}\ \varphi \subseteq
 K$ for each $\varphi \in B$ and with
 $$
 \sup_{(x,t)\in K} \{|\partial^\alpha\varphi(x,t)|:
 \varphi\in B,\, |\alpha|\leq 1\} <\infty \,.
 $$
 It suffices to show that there exist representatives $U_\eps$
 and $V_\eps$ of the solution (\ref{prvi sistem kolombo resenje})
 such that
 \be \label{prvi sistem asoc res1}
 \lim_{\eps\to 0}\sup_{\varphi\in B}|\int_{\R\times[0, \infty)}
 \Big((U_\eps)_t(x,t)+(U_\eps^2 - V_\eps)_x(x,t)\Big)\varphi(x,t)\,dx\,dt|=0,
 \ee
 and
 \be \label{prvi sistem asoc res2}
 \lim_{\eps\to 0}\sup_{\varphi\in B}|\int_{\R\times[0, \infty)}
 \Big((V_\eps)_t(x,t)+(\frac13 U_\eps^3- U_\eps)_x(x,t)\Big)\varphi(x,t)
 \,dx\,dt|=0.
 \ee
 Look first at (\ref{prvi sistem asoc res1}).
 Let $\varphi \in B$. Then
 \beast
 && \int_{\R\times [0,\infty)}
 \bigg[\Big((U_\eps(x,t)\Big)_t +
 \Big(U_\eps^2(x,t) - V_\eps(x,t)\Big)_x \bigg]
 \varphi(x,t)\,dx\,dt\\
 &&= \int_{\R\times [0,\infty)}
 \bigg[ \bigg\{
 G_\eps(x-ct)+s_{1}(t)
 (\alpha_{0}d^{-}(x-ct)+\alpha_{1}d^{+}(x-ct)) \bigg\}_t \\
 &&\quad
 + \bigg\{ G_\eps^2(x-ct)
 +2s_{1}(t)(\alpha_{0}u_0d^{-}(x-ct)+\alpha_{1}u_1d^{+}(x-ct)) \\
 &&\quad
 + s_{1}^2(t)(\alpha_{0}^2(d^{-})^2(x-ct)+\alpha_{1}^2(d^{+})^2(x-ct))
 -H_\eps(x-ct) \\
 &&\quad -s_{2}(t)(\beta_{0}D^{-}(x-ct)+\beta_{1}D^{+}(x-ct))
 \bigg\}_x
 \bigg] \varphi(x.t)\,dx\,dt\\
 && = \int_{\R\times [0,\infty)}
 \bigg[
 -c\pd_x
 G_\eps(x-ct)+s'_1(t)(\alpha_{0}d^{-}(x-ct)+\alpha_{1}d^{+}(x-ct))
 \\
 &&\quad - cs_1(t)\pd_x(\alpha_{0}d^{-}(x-ct)+\alpha_{1}d^{+}(x-ct))
 \\
 &&\quad +\pd_xG_\eps(x-ct) +
 s_{1}(t)\pd_x(\alpha_{0}u_0d^{-}(x-ct)+\alpha_{1}u_1d^{+}(x-ct))\\
 &&\quad + s_{1}^2(t)\pd_x(\alpha_{0}^2(d^{-})^2(x-ct)+
 \alpha_{1}^2(d^{+})^2(x-ct))\\
 &&\quad - \pd_xH_\eps(x-ct) - s_{2}(t)\pd_x(\beta_{0}D^{-}(x-ct)+
 \beta_{1}D^{+}(x-ct))
 \bigg]\varphi(x,t)\,dx\,dt\\
 \eeast
 \beast
 && = \int_{\R\times [0,\infty)}
 \underbrace{s'_1(t)(\alpha_{0}d^{-}(x-ct)+\alpha_{1}d^{+}(x-ct))}_{(1)}
 \varphi(x,t)\,dx\,dt \\
 &&\quad - \int_{\R\times [0,\infty)}
 \bigg[
 -\underbrace{cG_\eps(x-ct)}_{(2)} -
 \underbrace{cs_1(t)(\alpha_{0}d^{-}(x-ct)+\alpha_{1}d^{+}(x-ct))}_{(3)}
 \\
 &&\quad + \underbrace{G_\eps^2(x-ct)}_{(4)} + \underbrace{s_{1}(t)
 (\alpha_{0}u_0d^{-}(x-ct)+\alpha_{1}u_1d^{+}(x-ct))}_{(5)}\\
 &&\quad +\underbrace{s_{1}^2(t)(\alpha_{0}^2(d^{-})^2(x-ct)+
 \alpha_{1}^2(d^{+})^2(x-ct))}_{(6)}\\
 &&\quad - \underbrace{H_\eps(x-ct)}_{(7)}
 - \underbrace{s_{2}(t)(\beta_{0}D^{-}(x-ct)+
 \beta_{1}D^{+}(x-ct))}_{(8)}
 \bigg]
 \varphi_x(x,t)\,dx\,dt\\
 &&=(*)
 \eeast
 For the first member of this sum we have
 \beast
 && \int_{\R\times [0,\infty)}
 s'_1(t)(\alpha_{0}d^{-}(x-ct)+\alpha_{1}d^{+}(x-ct))\varphi(x,t)\,dx\,dt\\
 &&=\int_0^\infty \int_{\R}
 s'_1(t)\Big(\alpha_0\Big(-\frac{1}{2\eps}\phi\Big(
 \frac{-x+ct-4\varepsilon}{\eps}\Big) + \frac{1}{2\eps}\phi\Big(
 \frac{-x+ct-6\varepsilon}{\varepsilon} \Big) \Big)^{1/2} \\
 &&\quad +\alpha_1 \Big(\frac{1}{2\eps}\phi\Big(
 \frac{x-ct-4\varepsilon}{\eps}\Big) - \frac{1}{2\eps}\phi\Big(
 \frac{x-ct-6\varepsilon}{\varepsilon} \Big) \Big)^{1/2}
 \Big)\varphi(x,t)\,dx\,dt
 \eeast
 Here we used the $3'$SD-function from \cite{Ned}.
 Introducing suitable substitutions we obtain
 \beast
 &&\int_0^\infty s'_1(t)\int^{\infty}_{-1}
 \alpha_0\sqrt{\eps}\Big(\frac12\phi(z)\Big)^{1/2}
 \varphi(-\eps z -4\eps +ct, t)\,dz\,dt \\
 &&\quad -
 \int_0^\infty s'_1(t)\int^{\infty}_{-1}
 \alpha_0\sqrt{\eps}\Big(\frac12\phi(z)\Big)^{1/2}
 \varphi(-\eps z -6\eps +ct, t)\,dz\,dt \\
 && +\int_0^\infty s'_1(t)\int^{\infty}_{-1}
 \alpha_1\sqrt{\eps}\Big(\frac12\phi(z)\Big)^{1/2}
 \varphi(\eps z -4\eps +ct, t)\,dz\,dt \\
 &&\quad -
 \int_0^\infty s'_1(t)\int^{\infty}_{-1}
 \alpha_1\sqrt{\eps}\Big(\frac12\phi(z)\Big)^{1/2}
 \varphi(\eps z -6\eps +ct, t)\,dz\,dt.
 \eeast
 Applying the Lebesgue dominated convergence theorem two times
 successively to the corresponding sequences we conclude that
 this term tends to $0$ as $\eps \to 0$. A similar argument
 shows that each of the terms in the sum with the functions
 $d^\pm$ or
 $(d^\pm)^3$, i.e. (3) and (5), also goes to $0$ as
 $\eps \to 0$. So, look now at (2). We have
 \beast
 &&-\int_{\cR\times
 [0,\infty)}-cG_\eps(x-ct)\varphi_x(x,t)\,dx\,dt\\
 &&=\int_0^\infty \int_{-\infty}^{ct-\eps}cu_0
 \varphi_x(x,t)\,dx\,dt + \int_0^\infty \int_{ct+\eps}^\infty cu_1
 \varphi_x(x,t)\,dx\,dt\\
 &&=cu_0 \int_0^\infty \varphi(ct-\eps,t)\,dt -cu_1 \int_0^\infty
 \varphi(ct+\eps,t)\,dt\\
 &&\stackrel{\eps\to 0}{\longrightarrow}(cu_0-cu_1) \int_0^\infty
 \varphi(ct,t)\,dt\\
 &&=-c[G]\int_0^\infty \varphi(ct,t)\,dt,
 \eeast
 where we again applied the Lebesgue dominated convergence theorem.
 For (8) we obtain
 \beast
 &&-\int_{\cR\times [0,\infty)}
 -s_2(t)(\beta_{0}D^{-}(x-ct)+\beta_{1}D^{+}(x-ct))
 \varphi_x(x,t)\,dx\,dt\\
 &&=\int_0^\infty \int_{-\infty}^{ct-\eps}s_2(t)
 \frac{\beta_0}{\eps}\phi\Big(\frac{x-ct+2\eps}{\eps}\Big)
 \varphi_x(x,t)\,dx\,dt\\
 &&\quad +
 \int_0^\infty \int_{ct+\eps}^\infty s_2(t)
 \frac{\beta_1}{\eps}\phi\Big(\frac{x-ct-2\eps}{\eps}\Big)
 \varphi_x(x,t)\,dx\,dt\\
 &&=\int_0^\infty s_2(t)\int_{-\infty}^1 \beta_0\phi(z)
 \varphi_x(\eps z -2\eps +ct, t)\,dz\,dt \\
 &&\quad +
 \int_0^\infty s_2(t)\int^{\infty}_{-1}
 \beta_1\phi(z)\varphi_x(\eps z +2\eps +ct, t)\,dz\,dt\\
 &&\stackrel{\eps\to 0}{\longrightarrow}
 \int_0^\infty s_2(t) \beta_0
 \varphi_x(ct,t) \int_{-\infty}^1 \phi(z)\,dz\,dt
 +
 \int_0^\infty s_2(t)
 \beta_1 \varphi_x(ct,t) \int^{\infty}_{-1} \phi(z)\,dz\,dt\\
 && = (\beta_0+\beta_1)\int_0^\infty
 s_2(t)\varphi_x(ct,t)\,dt\\
 && = \int^\infty_0 s_2(t)\varphi_x(ct,t)\,dt.
 \eeast
 Repeating this for the remaining terms yields that (4)
 tends to $[G^2]\int_0^\infty \varphi(ct,t)\,dt$, (6) to
 $(\alpha_0^2+\alpha_1^2)\int_0^\infty s_1^2(t)\varphi_x(ct,t)\, dt$,
 and (7) to $-[H]\int_0^\infty \varphi(ct,t)\,dt$, as $\eps \to 0$.
 Therefore,
 \beast
 (*)\! &\! \stackrel{\eps\to 0}{\longrightarrow} \!&\!
   -c[G]\int_0^\infty \varphi(ct,t)\,dt
 +[G^2]\int_0^\infty \varphi(ct,t)\,dt\\
 \!&\!\!&\! \quad+ (\alpha_0^2+\alpha_1^2)\int^\infty_0 s_1^2(t)\varphi(ct,t)\,dt
 - [H] \int^\infty_0 \varphi(ct,t)\,dt
 - \int_0^\infty s_2(t)\varphi(ct,t)\,dt\\
 \!&\! =\! &\! 0,
 \eeast
 by (\ref{prvi sistem uslov1})-(\ref{prvi sistem
 uslov4}) and (\ref{prvi sistem asoc
 res1}) is satisfied. Similarly we conclude that it is also
 true for (\ref{prvi sistem asoc res2}). Hence, the solution
 (\ref{prvi sistem kolombo resenje}) is a $1$-strongly associated
 solution to (\ref{sistem1}).
 \ep\sms

 The projectable symmetry groups calculated at the beginning
 of this section transform $1$-strongly associated solutions to
 (\ref{sistem1}) to other $1$-strongly associated solutions,
 as shown by the following

 \bt \label{thm mm 1strogo asoc grup sim}
 The symmetry groups $G_1$ - $G_5$ of the system (\ref{sistem1}) are ${\cal
 AS}^1_\Delta$-symmetry groups.
 \et

 \pr
 By Theorem \ref{thm Q} it suffices to show that
 $G_1$ - $G_5$ are slowly increasing and have a factorization
 (\ref{faktorizacija matricna}) such that $Q$ depends only
 on $\eta$ and $x$. First we consider $G_1$. The action of $G_1$
 is given by
 \[
 g_\eta: (x,t,u,v) \to \left(e^{\eta}x,e^{\eta}t, u,v \right).
 \]
 Since $\Phi$ is the identity it follows that the map
 $$
 (u,v) \mapsto \Phi_g(x,t,u,v)
 $$
 is slowly increasing, uniformly for $x$ and $t$ in compact sets.
 It is easy to see that this is also true for the remaining
 groups. Next, for $G_1$ we have
 \beast
 \Delta_1(e^{-\eta}x,e^{-\eta}t, \prjedan u(e^{-\eta}x,e^{-\eta}t),
 \prjedan v(e^{-\eta}x,e^{-\eta}t))
 &\! =\! & e^{-\eta}u_t + 2e^{-\eta}uu_x - e^{-\eta}v_x \\
 &\! =\! & e^{-\eta}\Delta_1\\
 \Delta_2(e^{-\eta}x,e^{-\eta}t, \prjedan u(e^{-\eta}x,e^{-\eta}t),
 \prjedan v(e^{-\eta}x,e^{-\eta}t))
 &\! =\! & e^{-\eta}v_t +e^{-\eta}u^2u_x-e^{-\eta}u_x\\
 &\! =\! & e^{-\eta}\Delta_2,
 \eeast
 where $\Delta_1$ and $\Delta_2$ denote the first, respectively
 the second equation of the system (\ref{sistem1}).
 The matrix form of this factorization is given by
 $$
 \left[\begin{array}{c}
 \widetilde{\Delta}_1\\
 \widetilde{\Delta}_2
 \end{array}
 \right]=
 \left[\begin{array}{cc}
 e^{-\eta} & 0\\
 0 & e^{-\eta}
 \end{array}
 \right]\cdot
 \left[\begin{array}{c}
 \Delta_1\\
 \Delta_2
 \end{array}
 \right].
 $$
 Therefore, the matrix $Q$ depends only on $\eta$ and
 Theorem \ref{thm Q} provides that the $G_1$ is  ${\cal
 AS}^1_\Delta$-symmetry group. The factorizations for $G_2$-$G_5$ are
 $$
 G_2, G_4, G_5: \qquad
 \left[\begin{array}{c}
 \widetilde{\Delta}_1\\
 \widetilde{\Delta}_2
 \end{array}
 \right]=
 \left[\begin{array}{cc}
 1 & 0\\
 0 & 1
 \end{array}
 \right]\cdot
 \left[\begin{array}{c}
 \Delta_1\\
 \Delta_2
 \end{array}
 \right],
 $$
 $$
 G_3: \qquad\qquad\quad
 \left[\begin{array}{c}
 \widetilde{\Delta}_1\\
 \widetilde{\Delta}_2
 \end{array}
 \right]=
 \left[\begin{array}{cc}
 1 & 0\\
 \eta & 1
 \end{array}
 \right]\cdot
 \left[\begin{array}{c}
 \Delta_1\\
 \Delta_2
 \end{array}
 \right],
 $$
 hence again we conclude that these four symmetry groups
 are also ${\cal AS}^1_\Delta$-symmetry groups.
 \ep

\subsection{Zero Pressure Gas Dynamics Model}

 The next system of conservation laws we consider is given by

 \begin{equation} \begin{split}\label{model gasne dinamike}\displaystyle
 u_t + (uv)_x & = 0 \\
 (uv)_t + (uv^2)_x & = 0
 \end{split} \end{equation}

 with the same initial conditions (\ref{pocetni uslovi})
 as in the previous case. This Riemann problem is a zero pressure
 gas dynamics model, where $u$ is a density, hence nonnegative,
 and $v$ is a velocity.

 The quasilinear form of (\ref{model gasne dinamike}) is obtained by differentiating:

 \begin{equation} \begin{split}\displaystyle\label{model g d kvazilin}
 u_t + u_x v + uv_x & = 0 \\
 u_t v + uv_t + u_x v^2 + 2uvv_x & = 0.
 \end{split} \end{equation}

 If we compute $u_t$ from the first equation of
 (\ref{model g d kvazilin}) and insert into the second one we
 arrive to the following system:

 \begin{equation} \begin{split}\label{model g d redukovani}\displaystyle
 u_t + u_x v + uv_x & = 0 \\
 u(v_t + vv_x) & = 0.
 \end{split} \end{equation}

 From the second equation it can be seen that one possible solution is
 $u=0$, i.e. vacuum state. Therefore, we consider the other possibility
 $v_t+vv_x=0$, looking at the system

 \begin{equation} \begin{split}\label{siste2 posle delenja}
 u_t + u_x v + uv_x & = 0 \\
 v_t + vv_x & = 0.
 \end{split} \end{equation}

 The eigenvalues of this system are $\lambda_1(u,v) =
 \lambda_2(u,v) = v$, thus the system is weakly hyperbolic.
 The corresponding right eigenvector is $r(u,v)= [0 \; 1]^T$,
 so both characteristic fields are linearly degenerative.
 Hence, only contact discontinuities can appear as a solution
 and we calculate them: let $x=ct$ be a curve of discontinuities
 of the system (\ref{model g d redukovani}).
 Along this curve the Rankine-Hugoniot conditions must hold:
 \beast
 c[u]&=&[uv]\\
 c[uv]&=&[uv^2].
 \eeast
 Equating $c$ from these equations yields
 $$
 [u][uv^2]=[uv]^2,
 $$
 so
 $$
 u_0u_1(v_1-v_0)^2=0.
 $$
 Therefore, there exist three solutions: $u_0=0$, $u_1=0$ and $v_0=v_1$.
 For $u_0=0$ we calculate $c=v_1$, $\lambda_i(u_0,v_0)=v_0$
 and $\lambda_i(u_1,v_1)=v_1$, $i=1,2$ and similarly for
 $u_1=0$. In the third case $c=v_0=v_1$ and also $\lambda_i(u_0,v_0)=
\lambda_i(u_1,v_1)=v_0=v_1$, thus the Lax entropy condition
 (\ref{Laxov uslov}) is satisfied. Hence the initial conditions
 $(u_0,v_0)$ and $(u_1,v_1)$ can be connected by contact discontinuity
 only when $v_0=v_1$. Finally, combining contact discontinuities and
 vacuum states we obtain the classical solution of the Riemann
 problem (\ref{model gasne dinamike}) when $v_0<v_1$:
 $$
 (u,v)(x,t)=\left\{\begin{array}{ll}
            (u_0,v_0), & \quad \frac{x}{t}<v_0\\
            (0,v(x,t)), &\quad v_0\leq \frac{x}{t}\leq v_1\\
            (u_1,v_1), & \quad \frac{x}{t}>v_1.
            \end{array}
            \right.
 $$
 In the case when $v_0>v_1$ this solution is not uniquely defined
 and certain nonregularities appear, which is studied in detail
 in \cite{Ned}. In that case the solution of (\ref{model gasne dinamike})
 is again a singular shock wave

 \begin{equation} \begin{split}
 \label{model g d kolombo resenje}
 U(x,t) & =G(x-ct)+s_{1}(t)(\alpha_{0}D^{-}+\alpha_{1}D^{+})
 +s_{2}(t)(\beta_{0}d^{-}+\beta_{1}d^{+}) \\
 V(x,t) & = H(x-ct)+s_{3}(t)(\gamma_{0}d^{-}+\gamma_{1}d^{+}),
 \end{split} \end{equation}

 where $G$ and $H$ are generalized step functions, $D$ and $d$
 are compatible S$\delta$- and $3$SD-functions 
 (cf.\ \cite{Ned}, Def.\ 3 with Ex.\ (i)) and
 \bea
 &&s_{1}(t)=\sigma_{1}t,\; \sigma_{1}=c[G]-[GH],\;\sigma_{1}>0,
     \label{model2}\\
 &&\alpha_{0}=\frac{v_{1}-c}{v_{1}-v_{0}},\;
   \alpha_{1}=\frac{c-v_{0}}{v_{1}-v_{0}},\label{model3}\\
 &&\sigma_{1}(\alpha_{0}v_{0}+\alpha_{1}v_{1})=\sigma_{1}c=c[GH]-[GH^{2}],
     \label{model4}\\
 &&-s_{2}(t)s_{3}^{2}(t)=s_{1}(t),\label{model5}\\
 &&\alpha_{0}(v_{0}^{2}-cv_{0})+\alpha_{1}(v_{1}^{2}-cv_{1})=
   \beta_{0}\gamma_{0}^{2}+\beta_{1}\gamma_{1}^{2}.\label{model6}
 \eea

 This time the function $s_1(t)$ denotes the strength of the
 singular shock wave.

 As for the first system (\ref{sistem1}) we are going to show
 that this solution is also a $1$-strongly associated solution to
 (\ref{model gasne dinamike}).

 \bt \label{thm sistem2 1strogo asoc}
 The solution (\ref{model g d kolombo resenje})
 of the system (\ref{model gasne dinamike}) is a $1$-strongly associated
 solution to (\ref{model gasne dinamike}).
 \et

 \pr
 The proof is similar to that for the system (\ref{sistem1}),
 so we have to show that there exist representatives
 $U_\eps$ and $V_\eps$ of the solutions $U$ and $V$
 defined in (\ref{model g d kolombo resenje}), such that for arbitrary
 set $B\subseteq \cbesk_c (\R\times [0,\infty))$ bounded in $\cjedan_c
 (\R\times [0,\infty))$ the following holds:
 \be \label{asoc res1}
 \lim_{\eps\to 0}\sup_{\varphi\in B}|\int_{\R\times[0, \infty)}
 \Big((U_\eps)_t(x,t)+(U_\eps V_\eps)_x(x,t)\Big)\varphi(x,t)\,dx\,dt|=0,
 \ee
 and for the second equation:
 \be \label{asoc res2}
 \lim_{\eps\to 0}\sup_{\varphi\in B}|\int_{\R\times[0, \infty)}
 \Big((U_\eps V_\eps)_t(x,t)+(U_\eps V_\eps^2)_x(x,t)\Big)\varphi(x,t)
 \,dx\,dt|=0.
 \ee
 Let $B$ be a bounded subset of $\cc^1(\R\times
 [0,\infty))$. First we prove (\ref{asoc res1}).
 Let $\varphi \in B$. Then
 \beast
 && \int_{\R\times [0,\infty)}
 \bigg[\Big((U_\eps(x,t)\Big)_t +
 \Big(U_\eps(x,t)V_\eps(x,t)\Big)_x \bigg]
 \varphi(x,t)\,dx\,dt \qquad\qquad\quad \\
 &&= \int_{\R\times [0,\infty)}
 \bigg[ \bigg\{
 G_\eps(x-ct)+s_{1}(t)
 (\alpha_{0}D^{-}(x-ct)+\alpha_{1}D^{+}(x-ct))
 \\
 &&\quad
 +s_{2}(t)(\beta_{0}d^{-}(x-ct)+\beta_{1}d^{+}(x-ct))
 \bigg\}_t + \bigg\{
 G_\eps(x-ct)H_\eps(x-ct)
 \\
 &&\quad
 +s_{1}(t)s_3(t)(\alpha_{0}D^{-}(x\!-\!ct)\!+\!\alpha_{1}D^{+}(x\!-\!ct))
 (\gamma_{0}d^{-}(x\!-\!ct)\!+\!\gamma_{1}d^{+}(x\!-\!ct)) \\
 &&\quad
 + s_{2}(t)s_3(t)(\beta_{0}d^{-}(x\!-\!ct)\!+\!\beta_{1}d^{+}(x\!-\!ct))
 (\gamma_{0}d^{-}(x\!-\!ct)\!+\!\gamma_{1}d^{+}(x\!-\!ct))\\
 &&\quad
 s_{1}(t)(\alpha_{0}v_0D^{-}(x-ct)+\alpha_{1}v_1D^{+}(x-ct))\\
 &&\quad +s_{2}(t)(\beta_{0}v_0d^{-}(x-ct)+\beta_{1}v_1d^{+}(x-ct))
 \\
 &&\quad +
 s_{3}(t)(\gamma_{0}u_0d^{-}(x-ct)+\gamma_{1}u_1d^{+}(x-ct))
 \bigg\}_x
 \bigg] \varphi(x.t)\,dx\,dt\\
 && = \int_{\R\times [0,\infty)}
 \bigg[
 -c\pd_x
 G_\eps(x-ct)+s'_1(t)(\alpha_{0}D^{-}(x-ct)+\alpha_{1}D^{+}(x-ct))
 \\
 &&\quad - cs_1(t)\pd_x(\alpha_{0}D^{-}(x-ct)+\alpha_{1}D^{+}(x-ct))
 \\
 &&\quad + s'_2(t)(\beta_{0}d^{-}(x-ct)+\beta_{1}d^{+}(x-ct)) \\
 &&\quad - cs_2(t)\pd_x(\beta_{0}d^{-}(x-ct)+\beta_{1}d^{+}(x-ct)) \\
 &&\quad +\pd_x (G_\eps(x-ct)H_\eps(x-ct))\\
 &&\quad + s_{1}(t)s_3(t)\pd_x
 \Big((\alpha_{0}D^{-}(x-ct)+\alpha_{1}D^{+}(x-ct))
 (\gamma_{0}d^{-}(x-ct) \\
 &&\qquad
 +\gamma_{1}d^{+}(x-ct))\Big)\\
 &&\quad + s_{2}(t)s_3(t)\pd_x
 \Big((\beta_{0}d^{-}(x-ct)+ \beta_{1}d^{+}(x-ct))
 (\gamma_{0}d^{-}(x-ct)\\
 &&\qquad +\gamma_{1}d^{+}(x-ct))\Big)\\
 &&\quad+s_1(t)\pd_x(\alpha_{0}v_0D^{-}(x-ct)+\alpha_{1}v_1D^{+}(x-ct))
 \\
 &&\quad + s_{2}(t)\pd_x(\beta_{0}v_0d^{-}(x-ct)+\beta_{1}v_1d^{+}(x-ct))
 \\
 &&\quad +s_{3}(t)\pd_x(\gamma_{0}u_0d^{-}(x-ct)+\gamma_{1}u_1d^{+}(x-ct))
 \bigg]\varphi(x,t)\,dx\,dt\\
 && = \int_{\R\times [0,\infty)}
 \bigg[
 \underbrace{s'_1(t)(\alpha_{0}D^{-}(x-ct)+\alpha_{1}D^{+}(x-ct))}_{(1)}
 \\
 &&\quad
 \underbrace{s'_2(t)(\beta_{0}d^{-}(x-ct)+\beta_{1}d^{+}(x-ct))}_{(2)}
 \bigg] \varphi(x,t)\,dx\,dt \\
 \eeast
 \beast
 &&\quad - \int_{\R\times [0,\infty)}
 \bigg[
 -\underbrace{cG_\eps(x-ct)}_{(3)} -
 \underbrace{cs_1(t)(\alpha_{0}D^{-}(x-ct)+\alpha_{1}D^{+}(x-ct))}_{(4)}
 \\
 &&\quad - \underbrace{cs_2(t)(\beta_{0}d^{-}(x-ct)+\beta_{1}d^{+}(x-ct))}_{(5)}
 +\underbrace{G_\eps(x-ct)H_\eps(x-ct)}_{(6)}\\
 &&\quad +\underbrace{s_{1}(t)s_3(t)(\alpha_{0}D^{-}(x\!-\!ct)
 \!+\!\alpha_{1}D^{+}(x\!-\!ct))
 (\gamma_{0}d^{-}(x\!-\!ct)\!+\!\gamma_{1}d^{+}(x\!-\!ct))}_{(7)}\\
 && \quad +
 \underbrace{s_{2}(t)s_3(t)(\beta_{0}d^{-}(x\!-\!ct)
 \!+\!\beta_{1}d^{+}(x\!-\!ct))
 (\gamma_{0}d^{-}(x\!-\!ct)\!+\!\gamma_{1}d^{+}(x\!-\!ct))}_{(8)}\\
 &&\quad +\underbrace{s_1(t)(\alpha_{0}v_0D^{-}(x-ct)+\alpha_{1}v_1D^{+}(x-ct))}_{(9)}
 \\
 &&\quad +\underbrace{s_{2}(t)(\beta_{0}v_0d^{-}(x-ct)+\beta_{1}v_1d^{+}(x-ct))}_{(10)}
 \\
 &&\quad + \underbrace{s_{3}(t)(\gamma_{0}u_0d^{-}(x-ct)+
 \gamma_{1}u_1d^{+}(x-ct))}_{(11)}
 \bigg]
 \varphi_x(x,t)\,dx\,dt\\
 &&=(*)
 \eeast

 Consider now each of the terms in the last sum,
 like in the proof of Theorem
 \ref{thm sistem1 1strogo asoc}. For (1) we have:
 \beast
 && \int_{\R\times [0,\infty)}
 s'_1(t)(\alpha_{0}D^{-}(x-ct)+\alpha_{1}D^{+}(x-ct))\varphi(x,t)\,dx\,dt\\
 &&=\int_0^\infty \int_{\R}
 s'_1(t)\Big(\frac{\alpha_0}{\eps}\phi\Big(\frac{x-ct+2\eps}{\eps}\Big)+
 \frac{\alpha_1}{\eps}\phi\Big(\frac{x-ct-2\eps}{\eps}\Big)\Big)\varphi(x,t)\,dx\,dt
 \eeast
 Next, split up this integral as a sum of two integrals
 and take suitable substitutions. This yields
 \beast
 &&\int_0^\infty \int_{-\infty}^1
 s'_1(t)\alpha_0\phi(z)\varphi(\eps z -2\eps +ct, t)\,dz\,dt \\
 &&\quad +\int_0^\infty \int^{\infty}_{-1}
 s'_1(t)\alpha_1\phi(z)\varphi(\eps z +2\eps +ct, t)\,dz\,dt
 \eeast
 Here we apply the Lebesgue dominated convergence
 theorem, first to the sequences $\big(\int_{-\infty}^1 \varphi
 (\eps z -2\eps +ct, t)\,dz\big)_\eps$ and
 $\big(\int^{\infty}_{-1} \varphi(\eps z +2\eps +ct, t)
 \,dz\big)_\eps$, and then to the sequences $\big(\varphi(\eps
 z -2\eps +ct, t)\big)_\eps$ and $\big(\varphi(\eps z +2\eps +ct,
 t)\big)_\eps$. Then the last integral
 \beast
 && \stackrel{\eps\to 0}{\longrightarrow} \int_0^\infty s'_1(t) \alpha_0
 \varphi(ct,t) \int_{-\infty}^1 \phi(z)\,dz\,dt + \int_0^\infty s'_1(t)
 \alpha_1 \varphi(ct,t) \int^{\infty}_{-1} \phi(z)\,dz\,dt\\
 && = (\alpha_0+\alpha_1)\int_0^\infty
 s'_1(t)\varphi(ct,t)\,dt\\
 && = \int^\infty_0s'_1(t)\varphi(ct,t)\,dt,
 \eeast
 since by definition of S$\delta$-functions, $\int \phi(z)\,dz=1$
 on the domain of $\phi$ (and that is the interval $[-1,1]$) and
 $\alpha_0+\alpha_1 =1$ by assumption. It is obvious that a
 procedure similar to this one and those from the proof of
 Theorem \ref{thm sistem1 1strogo asoc} is repeated for each of
 the remaining $9$ terms (some of them are explicitly calculated in
 the proof of Theorem \ref{thm sistem1 1strogo asoc}).
 The only difference is that for the solution of the system
 (\ref{model g d kolombo resenje}) we used $3$SD-functions instead
 of $3'$SD-functions from the solution of (\ref{sistem1}).
 This implies that all terms with $d^{\pm}$ or $(d^{\pm})^2$ tend to
 $0$ as $\eps \to 0$.

 Joining all together we have
 \beast
 (*) \!&\! \stackrel{\eps\to 0}{\longrightarrow} \!&\!
 \int^\infty_0s'_1(t)\varphi(ct,t)\,dt -c[G]\int_0^\infty
 \varphi(ct,t)\,dt + \int^\infty_0 cs_1(t)\varphi_x(ct,t)\,dt \\
 \!&\!\!&\! \quad + [GH]\int_0^\infty \varphi(ct,t)\,dt - \int^\infty_0
 (\alpha_0v_0+\alpha_1v_1)s_1(t)\varphi_x(ct,t)\,dt\\
 \!&\! =\!&\! 0,
 \eeast
 since from (\ref{model2}) $s'_1(t)=c[G]-[GH]$ and
 $\alpha_0v_0+\alpha_1v_1=c$ from (\ref{model4}) and the
 condition $\alpha_0+\alpha_1=1$. Thus (\ref{model g d kolombo
 resenje}) is a $1$-strongly associated solution to the first equation
 of (\ref{model gasne dinamike}). Analogously, it can be seen that
 the same holds also in the case of the second equation,
 which proves the claim.
 \ep \sms

 The next task is to calculate symmetry groups of (\ref{model gasne
 dinamike}). Let us recall that (\ref{model gasne dinamike}),
 (\ref{model g d kvazilin}) and (\ref{model g d redukovani})
 are equivalent systems for all smooth solutions.
 Also, these systems are equivalent in the Colombeau algebra,
 since the elements of this algebra are equivalence classes of
 sequences of smooth functions. Therefore we look for the symmetry
 groups of the quasilinear system (\ref{model g d redukovani}).
 We start with the following

 \bt \label{thm sj}
 Let $\Delta(x,\un)=0$ and  $\Delta_i(x,\un)$, $i=1, \dots, k$,
 be nondegenerate differential equations on $M\subset X\times U$
 such that $\Delta$ can be written as a product
 \be \label{sj}
 \Delta = \prod_{i=1}^k \Delta_i.
 \ee
 If we denote the corresponding algebras of infinitesimal generators
 of symmetries of $\Delta$ and $\Delta_i$ by ${\goth g}$ and
 ${\goth g}_i$ respectively, $i=1, \dots, k$, then
 \be \label{podskup}
 \bigcap\limits_{i=1}^k {\goth g}_i \subseteq {\goth g}.
 \ee
 \et

 \pr
 Let $\vp \in \cap_{i=1}^k {\goth g}_i$. Then $\vp \in {\goth g}_i$, 
 $\forall i=1,\dots, k$, i.e. $\vp$ is a generator of a local 
 one-parameter symmetry group of each equation
 $$
 \Delta_i(x, \un)= 0, \qquad i=1,\dots, k.
 $$ 
 By the infinitesimal criterion and (\ref{ekv infkr}) 
 we may write
 \be \label{ffaktorizacija}
 \prn \vp(\Delta_i(x, \un))= Q_{i}\cdot
 \Delta_{i}(x, \un), \quad i=1, \dots,l,
 \ee
 with well-defined functions $Q_{i}$, $i=1,\dots,l$. 
 Since $\prn \vp$ is a vector field on the $n$-jet space $M^{(n)}$,
 the Leibniz rule for the product derivative yields
 \beast
 && \prn \vp (\Delta)=\prn \vp \big(\prod_{i=1}^k \Delta_i\big)=
 \sum\limits_{i=1}^k \Delta_1 \cdot \ldots \cdot\prn \vp
 (\Delta_i)\cdot \ldots \cdot \Delta_k \\
 && \qquad\qquad = \sum\limits_{i=1}^k \Delta_1 \cdot \ldots\cdot
 Q_i\Delta_i \cdot\ldots\cdot \Delta_k = Q\cdot \Delta, 
 \eeast
 where $Q= Q_1 + \ldots +Q_l$.
 Another application of Theorem \ref{thm infinitesimalni
 kriterijum} provides that $\vp \in {\goth g}$, which proves the claim.
 \ep \sms

 According to this theorem, the intersection of the symmetry groups of the
 system (\ref{siste2 posle delenja}) with
 the symmetry groups of $u=0$ will provide symmetry groups of the
 system (\ref{model g d redukovani}). (It should be noticed that due to the
 inclusion in (\ref{podskup}) not all symmetry groups of 
 (\ref{model g d redukovani}) will be obtained. However, Theorem \ref{thm sj}
 is of great help, since a direct computation of the symmetry groups
 of (\ref{model g d redukovani}) is a very difficult task.)
 The symmetry groups of $u=0$ can easily be calculated.
 Namely, by the infinitesimal criterion
 (\ref{infinitesimalni kriterijum}) it follows that the infinitesimal
 generators of $u=0$ are obtained as solutions of
 $\vp(u) = \phi = 0$, whenever $u=0$ (we assumed here that
 $\vp=\xi(x,t,u,v) \partial_x + \tau(x,t,u,v) \partial_t +
 \phi(x,t,u,v)\partial_u + \psi(x,t,u,v)\partial_v$).

 Now we follow the procedure for calculating symmetry groups
 of the system (\ref{siste2 posle delenja}).

 (1)
 \[
 {\bf w} = \xi(x,t,u,v) \partial_x + \tau(x,t,u,v) \partial_t +
 \phi(x,t,u,v)\partial_u + \psi(x,t,u,v)\partial_v.
 \]

 (2) The first prolongation is given by
 \[
 \mbox{\rm pr}^{(1)}{\bf w} = {\bf w} + \phi^x \partial_{u_x} +
 \phi^t \partial_{u_t} + \psi^x \partial_{v_x} + \psi^t
 \partial_{v_t},
 \]
 where $\phi^x, \phi^t, \psi^x$ and $\psi^t$ are the same as in (\ref{guzva}).

 (3) Since the equations of this systems are
 \beast
 & \Delta_1(x,t,u,v,u_x,v_x,u_t,v_t) = u_t + u_x v + u v_x &\\
 & \Delta_2(x,t,u,v,u_x,v_x,u_t,v_t) = v_t + v v_x, &
 \eeast
 we have to solve
 \beast
 & \mbox{\rm pr}^{(1)}{\bf w}(\Delta_1) = \phi^t + v \phi^x + \phi u_x
      + \phi v_x + u \psi^x = 0 &\\
 & \mbox{\rm pr}^{(1)}{\bf w}(\Delta_2) = \psi^t + \psi v_x +
       v \psi^x = 0. &
 \eeast
 Again we look only for the projectable symmetry groups.
 Inserting (\ref{guzva}), having in mind that the partial derivatives of
 $\xi$ and $\tau$ with respect to $u$ and $v$ vanish, and then
 substituting $u_t$ by $-u_x v - u v_x$ and $v_t$ by $-v v_x$,
 we obtain
 \beast
 && \phi_t - \phi_uvu_x - \phi_uuv_x - \phi_vvv_x - \xi_tu_x +
 \tau_tvu_x + \tau_tuv_x + \phi_xv  \\
 && +\, \phi_uvu_x + \phi_vvv_x - \xi_xvu_x + \tau_xv^2u_x +
 \tau_xuvv_x + \psi u_x  + \phi v_x + \psi_xu \\
 && +\, \psi_uuu_x+ \psi_vuv_x - \xi_xuv_x -\tau_xuvv_x = 0 \\
 && \psi_t - \psi_uvu_x - \psi_uuv_x - \psi_vvv_x -
 \xi_tv_x + \tau_tvv_x + \psi v_x + \psi_xv \\
 && +\, \psi_uvu_x + \psi_vvv_x - \xi_xvv_x + \tau_xv^2v_x = 0.
 \eeast
 (4) Coefficients of $1$, $u_x$ and $v_x$ equating with
 $0$ yield the following equations
 \bea
 && \phi_t + v\phi_x + u\psi_x = 0 \nonumber\\
 && -\xi_t + \tau_t v - \xi_x v + \tau_x v^2 + \psi +u\psi = 0 \nonumber\\
 && -\phi_uu + \tau_tu +\tau_xuv + \phi + \psi_vu - \xi_xu +\tau_xuv =0 \nonumber\\
 && \psi_t + \psi_xv = 0 \nonumber\\
 && 0 =0 \nonumber\\
 && -\psi_uu  - \xi_t + \tau_tv + \psi - \xi_xv + \tau_xuv = 0. \nonumber
 \eea
 (5) The solution is
 \bea
 && \xi(x,t) = c_1 + c_3 t + (c_2+c_8t)x+c_5x^2 \nonumber \\
 && \tau(x,t) = c_6 + (c_7+c_8t)t +(c_4+c_5t)x  \nonumber \\
 && \phi(x,t,u,v) = u\alpha (x,t,v)  \nonumber \\
 && \psi(x,t,u,v) = c_3 \! + \! c_2v \! + \! 2c_5xv \! + \!
 c_8(tv\!+\!x) \! - \! v(c_7\!+\!2c_8t\!+\!c_5x) \! - \!
     v^2(c_4\!+\!c_5t), \nonumber
 \eea
 where $\alpha$ is a function which depends on $x,t$ and $v$ and
 satisfies the equation
 \be \label{bilo 540}
 \frac{c_8+c_5v+\alpha_t+v\alpha_x}{v}=0.
 \ee
 The eight constants $c_1 - c_8$ generate eight linearly independent
 infinitesimal generators of one-parameter projectable symmetry groups,
 while $\alpha(x,t,v)$ generates an infinite-dimensional group. From
 (\ref{bilo 540}) we see that $\alpha(x,t,v)$ must depend on
 constants $c_5$ and $c_8$. It is also clear that the function
 $\alpha$ is not uniquely determined. Hence, in order to calculate the
 infinitesimal generators of the projectable symmetry groups we choose
 one possibility for $\alpha$:
 $$
 \alpha(x,t,v)= -c_5 x - c_8 t + \beta(v).
 $$

 Now we can write all infinitesimal generators:
 \beast
 && {\bf w}_1 = \partial_x  \\
 && {\bf w}_2 = x\partial_x + v\partial_v \\
 \eeast
 \beast
 && {\bf w}_3 = t\partial_x + \partial_v \\
 && {\bf w}_4 = x\partial_t - v^2\partial_v \\
 && {\bf w}_5 = x^2\partial_x+xt\partial_t-xu\partial_u
     + (xv-tv^2)\partial_v \\
 && {\bf w}_6 = \partial_t \\
 && {\bf w}_7 = t\partial_t - v\partial_v \\
 && {\bf w}_8 = xt\partial_x+t^2\partial_t - tu\partial_u
     + (x-tv)\partial_v \\
 && {\bf w}_{\beta} = u\beta(v)\partial_u.
 \eeast

 (6) The one-parameter transformation groups generated by the
 vector fields $\vpw_1 - \vpw_8$ and $\vpw_\beta$ are:

 \be \label{ggg}\displaystyle
 \begin{array}{cl}
  G_1: &  (x,t,u,v) \to \displaystyle \left(x+\eta, t, u, v\right) \\
  G_2: &  (x,t,u,v) \to \displaystyle \left(e^{\eta}x,t,
 u,e^{\eta} v \right) \\
  G_3: &  (x,t,u,v) \to \displaystyle \left(x+\eta t,t,
 u, v+\eta \right)\\
  G_4: &  (x,t,u,v) \to \displaystyle \left(x,t+\eta x,
 u,\frac{v}{1+\eta v}\right) \\
  G_5: &  (x,t,u,v) \to \displaystyle \left(\frac{x}{1-\eta x},\frac{t}{1-\eta x},
 (1-\eta x)u,\frac{v}{1-\eta (x-tv)}\right) \\
  G_6: &  (x,t,u,v) \to \displaystyle \left(x,t+\eta,
 u,v\right) \\
  G_7: &  (x,t,u,v) \to \displaystyle \left(x,e^{\eta} t,
 u,e^{-\eta} v\right) \\
  G_8: &  (x,t,u,v) \to \displaystyle \left(\frac{x}{1-\eta t},\frac{t}{1-\eta t},
 (1-\eta t)u,\eta x + (1-\eta t)v\right) \\
  G_\beta: & (x,t,u,v) \to
  \displaystyle \left(x,t,e^{\eta\beta(v)} u,v\right)
 \end{array}
 \ee
 Since each of the groups in (\ref{ggg}) is a symmetry group of
 the system (\ref{siste2 posle delenja}), from (\ref{ftilda projektabilna})
 it follows that if $u$ and $v$ are solutions so are the functions

 $$
 \begin{array}{lll}
 (1) & \widetilde{u}: (x,t) \to \displaystyle u(x -\eta, t) &
       \widetilde{v}: (x,t) \to \displaystyle v(x -\eta, t) \\
 (2) & \widetilde{u}: (x,t) \to \displaystyle u(e^{-\eta}x, t) &
       \widetilde{v}: (x,t) \to \displaystyle e^{\eta} v(e^{-\eta}x,t) \\
 (3) & \widetilde{u}: (x,t) \to \displaystyle u(x-\eta t, t) &
       \widetilde{v}: (x,t) \to \displaystyle v(x - \eta t, t) + \eta  \\
 (4) & \widetilde{u}: (x,t) \to \displaystyle u(x, t-\eta x) &
       \widetilde{v}: (x,t) \to \displaystyle \frac{v(x, t-\eta x)}{1+\eta v(x,
       t-\eta x)}\\
 (5) & \widetilde{u}: (x,t)\!\to\! \displaystyle \frac{u\left(\displaystyle
       \frac{x}{1+\eta x},
       \displaystyle \frac{t}{1+\eta x}\right)}{1+\eta x} &\!\!\!\!
       \widetilde{v}: (x,t)\! \to\! \displaystyle \frac{(1+\eta x)v\left(\displaystyle
       \frac{x}{1+\eta x},\displaystyle \frac{t}{1+\eta x}\right)}
       {1+\eta tv\left(\displaystyle \frac{x}{1+\eta x},\displaystyle \frac{t}{1+\eta x}\right)} \\
 (6) & \widetilde{u}: (x,t) \to \displaystyle u(x, t-\eta) &
       \widetilde{v}: (x,t) \to \displaystyle v(x, t-\eta) \\
 (7) & \widetilde{u}: (x,t) \to \displaystyle u(x, e^{-\eta}t) &
       \widetilde{v}: (x,t) \to \displaystyle e^{-\eta} v(x, e^{-\eta}t) \\
 (8) & \widetilde{u}: (x,t) \to \displaystyle \frac{u\left(\displaystyle \frac{x}{1+\eta t},
       \displaystyle \frac{t}{1+\eta t}\right)}{1+\eta t} &
       \widetilde{v}: (x,t) \to \displaystyle \frac{\eta x+v\left(\displaystyle \frac{x}{1+\eta t},
       \displaystyle \frac{t}{1+\eta t}\right)}{1+\eta t} \\
 (\beta) & \widetilde{u}: (x,t) \to \displaystyle e^{\eta\beta(v)}u(x,t) &
       \widetilde{v}: (x,t) \to v(x,t)
 \end{array}
 $$

 From the remark given after Theorem \ref{thm sj} it follows that all
 calculated symmetry groups are also symmetry groups of the system
 (\ref{model g d redukovani}): for infinitesimal generators ${\bf
 w}_5$, ${\bf w}_8$ and ${\bf w}_\beta$ the coefficients of $\pd_u$ are
 $\frac{1}{1+\eta x} u$, $\frac{1}{1+\eta t} u$ and $ e^{\eta\beta(v)}u$
 respectively, hence they vanish when $u=0$, while for the rest $\phi=0$.

 The matrix factorizations of (\ref{model gasne dinamike}) with respect to
 symmetry groups in (\ref{ggg}) are:

 $$
 G_1, G_6: \qquad
 \left[\begin{array}{c}
 \widetilde{\Delta}_1\\
 \widetilde{\Delta}_2
 \end{array}
 \right]=
 \left[\begin{array}{cc}
 1 & 0\\
 0 & 1
 \end{array}
 \right]\cdot
 \left[\begin{array}{c}
 \Delta_1\\
 \Delta_2
 \end{array}
 \right]
 $$
 $$
 G_2: \qquad\qquad
 \left[\begin{array}{c}
 \widetilde{\Delta}_1\\
 \widetilde{\Delta}_2
 \end{array}
 \right]=
 \left[\begin{array}{cc}
 1 & 0\\
 0 & e^\eta
 \end{array}
 \right]\cdot
 \left[\begin{array}{c}
 \Delta_1\\
 \Delta_2
 \end{array}
 \right]
 $$
 $$
 G_3: \qquad\qquad
 \left[\begin{array}{c}
 \widetilde{\Delta}_1\\
 \widetilde{\Delta}_2
 \end{array}
 \right]=
 \left[\begin{array}{cc}
 1 & 0\\
 \eta & 1
 \end{array}
 \right]\cdot
 \left[\begin{array}{c}
 \Delta_1\\
 \Delta_2
 \end{array}
 \right]
 $$
 $$
 G_4: \qquad
 \left[\begin{array}{c}
 \widetilde{\Delta}_1\\
 \widetilde{\Delta}_2
 \end{array}
 \right]=
 \left[\begin{array}{cc}
 \frac{1+2\eta v}{(1+\eta v)^2} &  \frac{-\eta}{(1+\eta v)^2} \\
 \frac{2\eta v^2}{(1+\eta v)^3} & \frac{1-\eta v}{(1+\eta v)^3}
 \end{array}
 \right]\cdot
 \left[\begin{array}{c}
 \Delta_1\\
 \Delta_2
 \end{array}
 \right]
 $$
 $$
 G_5: \qquad
 \left[\begin{array}{c}
 \widetilde{\Delta}_1\\
 \widetilde{\Delta}_2
 \end{array}
 \right]=
 \left[\begin{array}{cc}
 \frac{1+2\eta tv}{(1+\eta x)^2(1+\eta tv)} &  \frac{-\eta t}{(1+\eta x)^2(1+\eta tv)} \\
 \frac{2\eta tv^2}{(1+\eta x)(1+\eta tv)^3} & \frac{1-\eta tv}{(1+\eta x)(1+\eta tv)^3}
 \end{array}
 \right]\cdot
 \left[\begin{array}{c}
 \Delta_1\\
 \Delta_2
 \end{array}
 \right]
 $$
 $$
 G_7: \qquad\qquad
 \left[\begin{array}{c}
 \widetilde{\Delta}_1\\
 \widetilde{\Delta}_2
 \end{array}
 \right]=
 \left[\begin{array}{cc}
 e^{-\eta} & 0\\
 0 & e^{-2\eta}
 \end{array}
 \right]\cdot
 \left[\begin{array}{c}
 \Delta_1\\
 \Delta_2
 \end{array}
 \right]
 $$
 $$
 G_8: \qquad\qquad
 \left[\begin{array}{c}
 \widetilde{\Delta}_1\\
 \widetilde{\Delta}_2
 \end{array}
 \right]=
 \left[\begin{array}{cc}
 \frac{1}{(1+\eta t)^3} & 0 \\
 \frac{\eta x}{(1+\eta t)^4} & \frac{1}{(1+\eta t)^4}
 \end{array}
 \right]\cdot
 \left[\begin{array}{c}
 \Delta_1\\
 \Delta_2
 \end{array}
 \right]
 $$
 $$
 G_\beta:
 \left[\begin{array}{c}
 \widetilde{\Delta}_1\\
 \widetilde{\Delta}_2
 \end{array}
 \right]=
 \left[\begin{array}{cc}
 (1-\eta v\beta'(v))e^{\eta \beta(v)} & \eta \beta'(v)e^{\eta \beta(v)}\\
 -\eta uv^2\beta'(v)e^{\eta \beta(v)} & (1+\eta uv\beta'(v))e^{\eta \beta(v)}
 \end{array}
 \right]\cdot
 \left[\begin{array}{c}
 \Delta_1\\
 \Delta_2
 \end{array}
 \right].
 $$

 Therefore, the matrix of factorization $Q$ depends only on $x, t$ and $\eta$
 for all groups except for $G_4, G_5$ and $G_\beta$. In these three cases the
 factor $Q$ depends also on $v$. Beside that, from the transformed solutions
 we see that the groups $G_1$-$G_3$ and
 $G_6$-$G_8$ are slowly increasing, uniformly for $x$ and $t$ in compact sets.
 Thereby we have proved the next

 \bt
 The symmetry groups $G_1, G_2, G_3, G_6, G_7$ and $G_8$ of the system
 (\ref{model gasne dinamike}) are $1$-strongly associated symmetry groups,
 i.e. transform $1$-strongly associated solutions to 
 (\ref{model gasne dinamike}) to other $1$-strongly associated
 solutions.
 \et

 The remaining three groups $G_4, G_5$ and $G_\beta$ are not
 ${\cal AS}_\Delta^1$-groups for two reasons:
 first the condition that the map
 $$
 (u,v)\mapsto \Phi_g(x,t,u,v)
 $$
 is slowly increasing, uniformly for $x$ and $t$ in compact sets,
 does not hold globally. Second, the solution
 (\ref{model g d kolombo resenje})
 does not belong to the algebra $\cG_\infty$, which is necessary
 by Theorem \ref{thm Q}.

 Under certain assumptions on the solution $(u,v)$ defined in
 (\ref{model g d kolombo resenje}), and also on the parameter $\eta$, these
 problems can be avoided. Namely, if we assume that the function
 $v$ is nonnegative and $\eta \geq 0$ (then instead of a group we
 consider a semigroup) the symmetry groups $G_4$ and $G_5$ become
 slowly increasing, while for $G_8$ it should be supposed that
 $\beta(v)$ is a function of $L^\infty$-log-type.
 The second condition from Theorem \ref{thm Q} (i)
 would be fulfilled if we assume that the solution $(u,v)$ belongs
 to the algebra $\cG_\infty$.



 \end{document}